\documentclass[a4paper,12pt]{article}
\usepackage{amsmath,amsfonts,amssymb,amsthm}  
\newtheorem{Theorem}{Theorem}
\newtheorem{Lemma}{Lemma}
\newcommand{\ta}{\theta}
\newcommand{\la}{\lambda}

\author{Michel Lassalle\\
\small Centre National de la Recherche Scientifique\\[-0.8ex]
\small Institut Gaspard Monge, Universit\'e de Marne-la-Vall\'ee\\[-0.8ex]
\small 77454 Marne-la-Vall\'ee Cedex, France\\[-0.8ex]
\small \texttt{lassalle@univ-mlv.fr}\\[-0.8ex]
\small \texttt{http://igm.univ-mlv.fr/{\textasciitilde}lassalle}}

\title{A short proof of\\ generalized Jacobi-Trudi expansions\\
for Macdonald polynomials}

\begin{document}
\date{}
\maketitle

\begin{abstract}
We give an elementary proof of the development of Macdonald
polynomials in terms of ``modified complete'' and 
elementary symmetric functions.
\end{abstract}

\section{Introduction}

In spite of many efforts, the problem of finding an explicit 
analytic expansion for the ``zonal polynomials'' of Hua \cite{H} 
and James \cite{Ja} has kept unsolved for fifty years. 
This problem is now better understood in the 
more general framework of Macdonald
polynomials \cite{Ma}. Zonal polynomials are indeed a special case of
Jack polynomials, which are themselves a particular limit of
Macdonald polynomials.

Macdonald polynomials form a basis of the algebra
of symmetric functions with rational coefficients in two
parameters $q,t$. They generalize several classical bases of 
this algebra, including monomial, elementary, Schur, 
Hall--Littlewood, and Jack symmetric functions. These 
particular cases correspond to various specializations of the 
indeterminates $q$ and $t$. 

Although two combinatorial expansions are known for Macdonald 
polynomials, in terms of tableaux \cite{Ma} or 
determinants \cite{LLM1}, no analytic expansion was available 
for them. Such a development had only been worked out when the 
indexing partition is a hook~\cite{Ke1}, has length
two~\cite{JJ} or three~\cite{La},
and in the dual cases corresponding to parts
at most equal to $3$.

In \cite{LS1} a general solution has been announced for 
this problem, providing two explicit analytic developments
for Macdonald polynomials, in terms of elementary and
``modified complete'' symmetric functions. In the special case $q=t$, 
these two developments coincide 
with the classical Jacobi-Trudi formulas for Schur functions. 
Therefore they appear as generalized 
Jacobi-Trudi expansions for Macdonald polynomials.

In \cite{LS2} a proof was given, which relies on the Pieri formula 
for Macdonald polynomials. This Pieri formula defines an infinite 
multidimensional matrix, which can be explicitly inverted 
using a method originally developed by
Krattenthaler \cite{Kr1,Kr2}, and further adapted to  
the multivariate case by Schlosser \cite{Sch}.  

This inversion of the ``Pieri matrix'' was the natural problem to 
solve, and it was needed for bringing to light the explicit form of the 
solution. However this method is rather technical and essentially 
external to Macdonald's theory. 

The purpose of this paper is to give a completely different proof, 
which is very quick and elementary, does not use the 
general Pieri formula, and keeps totally inside the Macdonald framework. 

A strange fact is that this proof studies Macdonald polynomials 
$P_{\la}$, with $\la$ a partition of length $n$, by applying 
Macdonald theory to some auxiliary space of $n$ indeterminates 
$u_i, 1\le i \le n$, specialized to $u_i=q^{\la_i-\la_n}t^{n-i-1}$. 
It is an open question whether an analogous method can be 
applied to Macdonald polynomials 
associated with other root systems than $A_{n}$.

\section{Macdonald polynomials}

The standard reference for Macdonald polynomials is Chapter 6
of \cite{Ma}.

A partition $\la= (\la_1,...,\la_n)$ is a finite weakly decreasing
sequence of nonnegative integers, called parts. The number
$l(\la)$ of positive parts is called the length of
$\la$, and $|\la| = \sum_{i = 1}^{n} \la_i$
the weight of $\la$. For any integer $i\geq1$,
$m_i(\la) = \textrm{card} \{j: \la_j  = i\}$
is the multiplicity of the part $i$ in $\la$.  Clearly
$l(\la)=\sum_{i\ge1} m_i(\la)$ and
$|\la|=\sum_{i\ge1} im_i(\la)$. We shall also write
$\la= (1^{m_1},2^{m_2},3^{m_3},\ldots)$. We set
\[z_\la  = \prod_{i \ge  1} i^{m_i(\lambda)} m_i(\lambda) ! .\]
We denote $\la^{'}$ the partition conjugate to $\la$, whose parts
are given by $m_i(\la^{'})=\la_i-\la_{i+1}$. We
have $\la^{'}_i=\sum_{j \ge i} m_j(\la)$.

Let $X=\{x_1,x_2,x_3,\ldots\}$ be a (finite or infinite) set
of indeterminates and $\mathcal{S}$ the
corresponding algebra of symmetric functions
with coefficients in $\mathsf{Q}$. 
When $X$ is infinite,
elementary symmetric functions $e_{k}(X)$,
complete symmetric functions $h_{k}(X)$, 
power sum symmetric functions $p_{k}(X)$
form three algebraic bases of $\mathcal{S}$, which can thus be 
viewed as an abstract algebra over $\mathsf{Q}$ generated 
by the functions $e_{k}$, $h_{k}$ or $p_{k}$.

Let $q,t$ be two indeterminates. We define ${(a;q)}_{0}=1$ 
and for $k\ge 1$,
\[{(a;q)}_{k}=\prod_{i=0}^{k-1}(1-aq^i),
\qquad {(a;q)}_{\infty}=\prod_{i\geq0}(1-aq^i).\]
Let $\mathsf{Q}[q,t]$ be the field of rational functions  in
$q,t$, and $\mathsf{Sym}=\mathcal{S}\otimes\mathsf{Q}[q,t]$
the algebra of symmetric functions
with coefficients in $\mathsf{Q}[q,t]$. 

For any $k\ge 0$, the ``modified complete'' symmetric function
$g_{k}(X;q,t)$ is defined by the generating series
\begin{equation*}
\prod_{i \ge 1}
\frac{{(tux_i;q)}_{\infty}}{{(ux_i;q)}_{\infty}}
=\sum_{k\ge0} u^k g_{k}(X;q,t).
\end{equation*}
It is often written in $\la$-ring notation 
\cite[p.~223]{La1}, i.e.
\[g_{k}(X;q,t)=h_{k}\left[\frac{1-t}{1-q}\,X\right].\]

The symmetric functions $g_{k}(q,t)$ form an algebraic basis of
$\mathsf{Sym}$, which may be expanded in terms of any classical
basis. This development is explicitly given in \cite[pp.~311 and 314]{Ma} 
in terms of power sums and monomial symmetric functions,
and in \cite[Sect.~10, p.~237]{La1} in terms of other classical 
bases. 

We now restrict to the case of a finite set of 
indeterminates $X=\{x_1,\ldots,x_n\}$.
Let $T_{q,x_i}$ denote the $q$-deformation operator defined 
by
\[T_{q,x_i}f (x_1,\ldots,x_n)= f(x_1,\ldots,qx_i,\ldots,x_n),\]
and for any $1\le i \le n$, 
\[A_i(X;t) =\prod_{\begin{subarray}{c}k=1\\k\neq i\end{subarray}}^n
\frac{tx_i-x_k}{x_i-x_k}.\]

Macdonald polynomials $P_{\la}(X;q,t)$, with $\la$ a 
partition such that $l(\la) \le n$, are defined as the  
eigenvectors of the following difference operator
\[E(q,t)= \sum_{i=1}^{n} \, A_i(X;t) \,T_{q,x_i}.\]
One has
\[E(q,t) \, P_{\la}(X;q,t)=  \left(\sum_{i=1}^n q^{\la_i}\, 
t^{n-i}\right) P_{\la}(X;q,t).\]
More generally let $\Delta(X)=
\prod_{1\le i < j \le n} (x_i-x_j)$ be the Vandermonde determinant 
and $a$ some indeterminate. Macdonald polynomials $P_{\la}(X;q,t)$ 
are eigenvectors of the difference operator
\[D(a;q,t)=\frac{1}{\Delta(X)} \operatorname{det}_{1\le i,j \le n}
{\left[x_i^{n-j}\left(1+at^{n-j}T_{q,x_i}\right)\right]}.\]
One has
\[D(a;q,t) \, P_{\la}(X;q,t)=  \prod_{i=1}^n \Big(1+a\ q^{\la_i}\, 
t^{n-i}\Big) P_{\la}(X;q,t).\]

The polynomials
$P_{\la}(X;q,t)$ form an orthogonal basis of $\mathsf{Sym}$ 
with respect to the scalar product $<\, , \,>_{q,t}$ defined 
by
\[<p_\la,p_\mu>_{q,t}=\delta_{\la \mu}\,z_{\la}\,
\prod_{i=1}^{l(\la)} \frac{1-q^{\la_i}}{1-t^{\la_i}},\]
where $p_\mu$ denotes the power sum symmetric function 
$p_{\mu} = \prod_{k=1}^{l(\mu)} p_{\mu_k}$.
Let $Q_{\la}(X;q,t)$ denote the dual basis of $P_{\la}(X;q,t)$ for this 
scalar product.  

Equivalently if $Y=\{y_1,\ldots,y_n\}$ is another set of $n$ 
indeterminates, and
\[\Pi(X,Y;q,t)=
\prod_{i,j=1}^n\frac{{(tx_iy_j;q)}_{\infty}}{{(x_iy_j;q)}_{\infty}},\]
we have
\[\Pi(X,Y;q,t)=\sum_{\la}P_{\la}(X;q,t)Q_{\la}(Y;q,t).\]
This yields immediately
\begin{equation}
D(a;q,t)_{(X)}\,\Pi(X,Y;q,t)=D(a;q,t)_{(Y)}\, \Pi(X,Y;q,t),
\end{equation}
\begin{equation}
E(q,t)_{(X)}\, D(a;q,t)_{(X)}\, \Pi(X,Y;q,t)=
E(q,t)_{(Y)}\, D(a;q,t)_{(Y)}\, \Pi(X,Y;q,t),
\end{equation}
where the suffix $X$ (resp. $Y$) indicates operation on 
the $X$ (resp. $Y$) variables.
Now one has $P_{\la}(1/q,1/t)=P_{\la}(q,t)$ and 
$Q_{\la}(1/q,1/t)=(qt^{-1})^{|\la|}\ Q_{\la}(q,t)$, which 
implies
\begin{equation}
E(1/q,1/t)_{(X)}\, D(a;q,t)_{(X)}\, \Pi(X,Y;q,t)=
E(1/q,1/t)_{(Y)}\, D(a;q,t)_{(Y)}\, \Pi(X,Y;q,t).
\end{equation}

The Macdonald polynomials associated with a row or a column 
partition are given by  
\[P_{1^k}(q,t)=e_{k},\qquad Q_{(k)}(q,t)= g_{k}(q,t).\]
In the following, when the parameters $q,t$ 
are omitted, $P_{\mu}$
or $Q_{\mu}$ stands for $P_{\mu}(q,t)$ or $Q_{\mu}(q,t)$.

\section{Statement of results}

Let $\mathsf{N}$ be the set of
nonnegative integers. Let $u=(u_1,\ldots,u_n)$ be $n$ indeterminates 
and $\ta =(\ta_1,\ldots,\ta_n)\in \mathsf{N}^{n}$. For clarity
of notations, we
introduce $n$ auxiliary variables $v=(v_1,\ldots,v_n)$
defined by $v_k=q^{\ta_k}u_k$. We write
\begin{multline*}
C^{(q,t)}_{\ta_1,\ldots,\ta_n} (u_1,\ldots,u_n)= \prod_{k=1}^n
t^{\ta_k} \,\frac{(q/t;q)_{\ta_k}}{(q;q)_{\ta_k}}\,
\frac{(qu_k;q)_{\ta_k}}{(qtu_k;q)_{\ta_k}}\,
\prod_{1\le i < j \le n}
\frac{{(qu_i/tu_j;q)}_{\ta_i}}{{(qu_i/u_j;q)}_{\ta_i}} \,
\frac{{(tu_i/v_j;q)}_{\ta_i}}{{(u_i/v_j;q)}_{\ta_i}}\\
\times\frac{1}{\Delta(v)} \,
\operatorname{det}_{1\le i,j \le n}
\left[v_i^{n-j}
\left(1-t^{j-1} \frac{1-tv_i}{1-v_i}
\prod_{k=1}^n \frac{u_k-v_i}{tu_k-v_i}\right)\right].
\end{multline*}
Setting $u_{n+1}=1/t$ we have
\begin{multline*}
C^{(q,t)}_{\ta_1,\ldots,\ta_n} (u_1,\ldots,u_n)=
\prod_{1\le i< j\le n+1}
\frac{(qu_i/tu_j;q)_{\ta_i}}{(qu_i/u_j;q)_{\ta_i}}
\prod_{1\le i\le j\le n}
\frac{(tu_i/v_j;q)_{\ta_i}}{(u_i/v_j;q)_{\ta_i}}\\\times
\frac{1}{\Delta(v)} \,
\operatorname{det}_{1\le i,j \le n}
\left[v_i^{n-j}
\left(1-t^{j}
\prod_{k=1}^{n+1} \frac{u_k-v_i}{tu_k-v_i}\right)\right].
\end{multline*}

The following results were announced in~\cite{LS1} and 
proved in~\cite{LS2}.
\begin{Theorem}\label{theomain}
Let  $\la=(\la_1,...,\la_{n+1})$ be an arbitrary partition with length
$n+1$. For any $1\le k \le n+1$ define
$u_k=q^{\la_k-\la_{n+1}}t^{n-k}$. We have
\begin{equation*}
Q_{(\la_1,\ldots,\la_{n+1})}= \sum_{\ta\in\mathsf{N}^n}
C^{(q,t)}_{\ta_1,\ldots,\ta_n} (u_1,\ldots,u_n)\:
Q_{(\la_{n+1}-|\ta|)} \:
Q_{(\la_1+\ta_1,\ldots,\la_n+\ta_n)}.
\end{equation*}
\end{Theorem}

There exists an automorphism
$\omega_{q,t}={\omega_{t,q}}^{-1}$ of $\mathsf{Sym}$ such that
\[\omega_{q,t}(Q_{\la}(q,t))=P_{\la^{'}}(t,q),\qquad
\omega_{q,t}(g_{k}(q,t))=e_{k}.\]
Applying $\omega_{q,t}$ to Theorem 1, we obtain
the following equivalent result.
\begin{Theorem}\label{theodual}
Let  $\la=(1^{m_1},2^{m_2},\ldots,(n+1)^{m_{n+1}})$ be an arbitrary
partition consisting of parts at most equal to $n+1$.
For any $1\le k \le n+1$ define
$u_k=q^{n-k}t^{\sum_{j=k}^n m_j}$. We have
\begin{multline*}
P_{(1^{m_1},2^{m_2},\ldots,(n+1)^{m_{n+1}})}=
\sum_{\ta\in\mathsf{N}^n}
C^{(t,q)}_{\ta_1,\ldots,\ta_n} (u_1,\ldots,u_n) \:
e_{m_{n+1}-|\ta|} \\
\times 
P_{(1^{m_1+\ta_1-\ta_2},2^{m_2+\ta_2-\ta_3},\ldots,
{(n-1)}^{m_{n-1}+\ta_{n-1}-\ta_n},n^{m_n+m_{n+1}+\ta_n})}.
\end{multline*}
\end{Theorem}

The specialization $q=t$ corresponds to the case of Schur
functions. Then we have $P_\la(t,t)=Q_\la(t,t)=s_\la$ and
$g_{k}(t,t)=h_{k}$. It may be shown~\cite[Sect.~7]{LS2} that 
Theorem \ref{theomain} reads
\begin{equation*}
s_{(\la_1,\ldots,\la_{n+1})}= \sum_{\ta\in{(0,1)}^n}
(-1)^{|\ta|} \: h_{\la_{n+1}-|\ta|} \:
s_{(\la_1+\ta_1,\ldots,\la_n+\ta_n)}.
\end{equation*}
This is exactly the development of the Jacobi-Trudi 
determinant
\[s_\la=\operatorname{det}_{1\le i,j \le n+1} \,
[h_{\la_{i}-i+j}]\]
along its last row \cite[p.~41, (3.4)]{Ma}.
Conversely Theorem \ref{theodual} is exactly the development of 
the N\"{a}gelsbach-Kotska determinant
\[s_\la=\operatorname{det}_{1\le i,j \le n+1} \,
[e_{\la^{'}_{i}-i+j}]\]
along its last row. Our results thus appear as generalized 
Jacobi-Trudi expansions for Macdonald polynomials.

Let $\mathsf{M}^{(n)}$ denote the set 
of upper triangular $n \times n$ 
matrices with nonnegative integers, and $0$ on the diagonal.
By a straightforward iteration of Theorem \ref{theomain} 
we deduce immediately the analytic development of Macdonald 
polynomials in terms of the symmetric functions $g_k$.

\begin{Theorem}\label{anexpg}
Let  $\la=(\la_1,...,\la_{n+1})$ be an arbitrary partition with length
$n+1$. We have
\begin{multline*}
Q_{\la}(q,t)=
\sum_{\ta\in \mathsf{M}^{(n+1)}}
\prod_{k=1}^{n}
C_{\ta_{1,k+1},\ldots,\ta_{k,k+1}}^{(q,t)}
(\{u_i=q^{\la_i-\la_{k+1}+\sum_{j=k+2}^{n+1} 
(\ta_{i,j}-\ta_{k+1,j})}t^{k-i};1\le i\le k\}) \\ 
\times \prod_{k=1}^{n+1} g_{\la_{k}+\sum_{j=k+1}^{n+1}
\ta_{kj}-\sum_{j=1}^{k-1} \ta_{jk}}.
\end{multline*}
\end{Theorem}

This result may also be stated in terms of ``raising 
operators'' \cite[p.~9]{Ma}. For each pair of integers $1\le i<j \le n+1$ 
define an operator $R_{ij}$ acting on multi-integers $a=(a_1,\ldots,a_{n+1})$ by 
$R_{ij}(a)= (a_1,\ldots,a_i+1,\ldots,a_j-1,\ldots,a_{n+1})$. 
Any product $R= \prod_{i<j} R_{ij}^{\ta_{ij}}$, with 
$\ta=(\ta_{ij})_{1\le i<j \le n+1} \in \mathsf{M}^{(n+1)}$ is 
called a raising operator.
Its action may be extended to
any function $g_{\mu} = \prod_{k=1}^{n+1} g_{\mu_k}$, with 
$\mu$ a partition of length $n+1$, by setting $Rg_{\mu}= g_{R(\mu)}$.
In particular $R_{ij}g_{\mu}= g_{\mu_{1}}\ldots 
g_{\mu_{i}+1}\ldots g_{\mu_{j}-1}\ldots g_{\mu_{n+1}}$. 
Then the quantity appearing above in the right-hand side may 
be written
\[\prod_{k=1}^{n+1} g_{\la_{k}+\sum_{j=k+1}^{n+1}
\ta_{kj}-\sum_{j=1}^{k-1} \ta_{jk}}=
\left(\prod_{1\le i<j \le n+1}R_{ij}^{\ta_{ij}}\right) \,
g_{\la}.\]
For $q=t$ we recover the 
following variant of the Jacobi-Trudi expansion \cite[p.~42]{Ma}
\[s_{\la}= 
\left(\prod_{1\le i<j \le {n+1}} (1-R_{ij})\right) \,
h_{\la}.\]

Applying $\omega_{q,t}$, we immediatly deduce the following analytic expansion of
Macdonald polynomials in terms of elementary symmetric 
functions $e_k$.
\begin{Theorem}\label{anexpe}
Let  $\la=(1^{m_1},2^{m_2},\ldots,(n+1)^{m_{n+1}})$ be an arbitrary
partition consisting of parts at most equal to $n+1$.
We have
\begin{multline*}
P_{\la}(q,t)=  
\sum_{\ta\in \mathsf{M}^{(n+1)}}
\prod_{k=1}^n
C_{\ta_{1,k+1},\ldots,\ta_{k,k+1}}^{(t,q)}
(\{u_i=q^{k-i}t^{\sum_{j=i}^k m_j+\sum_{j=k+2}^{n+1} 
(\ta_{i,j}-\ta_{k+1,j})};1\le i\le k\}) \\  
\times \prod_{k=1}^{n+1}
e_{\,\sum_{j=k}^{n+1}m_j+\sum_{j=k+1}^{n+1} 
\ta_{kj}-\sum_{j=1}^{k-1}\ta_{jk}}.
\end{multline*}
\end{Theorem}
 
It is clear that the analytic developments given above are fully explicit.
It is worth considering these results in some particular 
cases \cite[p.~324]{Ma}, namely
$q=t$ (Schur functions), $q=1$ (elementary symmetric 
functions), $t=1$ (monomial symmetric functions),
$q=0$ (Hall--Littlewood  symmetric functions), and 
$q=t^\alpha, \ t\rightarrow1$ 
(Jack symmetric functions). This is done in detail in 
\cite{LS2}, together with the case of hook partitions 
$(r,1^s)$, already worked out by Kerov \cite{Ke1}.

\section{Proof} 

We shall only need the following two elementary lemmas. Let  $a$,
$u=(u_1,\ldots,u_{n+1})$ and $v=(v_1,\ldots,v_{n+1})$ be 
$2n+3$ indeterminates. Define
\[H(u;v)=\frac{1}{\Delta(v)} \,
\operatorname{det}_{1\le i,j \le n+1}
{\left[v_i^{n-j+1}
\left(1+at^{j}
\prod_{k=1}^{n+1} \frac{u_k-v_i}{tu_k-v_i}\right)\right]}.\]
With $v_i=q^{\ta_i}u_i, 1\le i\le n$, the last factor of $C^{(q,t)}_{\ta_1,\ldots,\ta_n} (u_1,\ldots,u_n)$ 
can be written
\[\frac{1}{\Delta(v)} \,
\operatorname{det}_{1\le i,j \le n}
{\left[v_i^{n-j}
\left(1-t^{j}
\prod_{k=1}^{n+1} \frac{u_k-v_i}{tu_k-v_i}\right)\right]}
=\lim_{a\rightarrow-1}\frac{1}{1+a} H(u,1/t;v,0).\]

\begin{Lemma}
$H(u,v)$ satisfies the two following functional equations :
\begin{enumerate}
\item[(i)]
\begin{multline*}
\sum_{i=1}^{n+1} 
\prod_{\begin{subarray}{c}k=1\\ k\neq i\end{subarray}}^{n+1}
\frac{v_i/t-v_k}{v_i-v_k}\,
\prod_{k=1}^{n+1}
\frac{1-v_i/u_k}{1-v_i/tu_k}
\, H(u;v_1,\ldots,v_i/q,\ldots,v_{n+1})=\\
\sum_{i=1}^{n+1} 
\prod_{\begin{subarray}{c}k=1\\ k\neq i\end{subarray}}^{n+1}
\frac{u_k/t-u_i}{u_k-u_i}\,
\prod_{k=1}^{n+1}
\frac{1-v_k/u_i}{1-v_k/tu_i}
\, H(u_1,\ldots,qu_i,\ldots,u_{n+1};v),
\end{multline*}
\item[(ii)]
\begin{multline*}
\sum_{i=1}^{n+1}  
\prod_{\begin{subarray}{c}k=1\\ k\neq i\end{subarray}}^{n+1}
\frac{tv_i-v_k}{v_i-v_k}\,
\prod_{k=1}^{n+1} \frac{1-qv_i/tu_k}{1-qv_i/u_k}
\, H(u;v_1,\ldots,qv_i,\ldots,v_{n+1})=\\
\sum_{i=1}^{n+1}  
\prod_{\begin{subarray}{c}k=1\\ k\neq i\end{subarray}}^{n+1}
\frac{tu_k-u_i}{u_k-u_i}\,
\prod_{k=1}^{n+1} \frac{1-qv_k/tu_i}{1-qv_k/u_i}
\, H(u_1,\ldots,u_i/q,\ldots,u_{n+1};v).
\end{multline*}
\end{enumerate}
\end{Lemma}
\begin{proof}
Consider the indeterminates
$X=(x_1,\ldots,x_{n+1})$ and $Y=(y_1,\ldots,y_{n+1})$, with 
$x_i=v_i$ and $y_i=1/u_i$. 
Writing $\Pi$ for $\Pi(X,Y;1/q,1/t)$ we have obviously
\begin{equation*}
\begin{split}
\Pi^{-1}\, T_{1/q,x_i}\, \Pi= \prod_{k=1}^{n+1} 
\frac{1-v_i/u_k}{1-v_i/tu_k}, &\qquad
\Pi^{-1}\, T_{1/q,y_i}\, \Pi= \prod_{k=1}^{n+1}
\frac{1-v_k/u_i}{1-v_k/tu_i},\\
H(u;v)=\Pi^{-1}\, D(a;1/q,1/t)&_{(X)}\, \Pi
=\Pi^{-1}\, D(a;1/q,1/t)_{(Y)}\, \Pi,
\end{split}
\end{equation*}
where the second equality follows from (1). 
Now we have
\begin{equation*}
\begin{split}
\Pi^{-1}\, T_{1/q,x_i}\, D(a;1/q,1/t)_{(X)}\, \Pi&=
\prod_{k=1}^{n+1} \frac{1-v_i/u_k}{1-v_i/tu_k}
\, T_{1/q,x_i} H(u;v)\\
&=\prod_{k=1}^{n+1} \frac{1-v_i/u_k}{1-v_i/tu_k}
\, H(u;v_1,\ldots,v_i/q,\ldots,v_{n+1}),
\end{split}
\end{equation*}
and similarly
\begin{equation*}
\Pi^{-1}\, T_{1/q,y_i}\, D(a;1/q,1/t)_{(Y)}\, \Pi=
\prod_{k=1}^{n+1} \frac{1-v_k/u_i}{1-v_k/tu_i}
\, H(u_1,\ldots,qu_i,\ldots,u_{n+1};v).
\end{equation*}
Thus (i) follows from (2). We have also
\begin{equation*}
\Pi^{-1}\, T_{q,x_i}\, \Pi= \prod_{k=1}^{n+1} 
\frac{1-qv_i/tu_k}{1-qv_i/u_k} \quad,\quad
\Pi^{-1}\, T_{q,y_i}\, \Pi= \prod_{k=1}^{n+1}
\frac{1-qv_k/tu_i}{1-qv_k/u_i},
\end{equation*}
which implies
\begin{equation*}
\begin{split}
\Pi^{-1}\, T_{q,x_i}\, D(a;1/q,1/t)_{(X)}\, \Pi&=
\prod_{k=1}^{n+1} \frac{1-qv_i/tu_k}{1-qv_i/u_k}
\, H(u;v_1,\ldots,qv_i,\ldots,v_{n+1}),\\
\Pi^{-1}\, T_{q,y_i}\, D(a;1/q,1/t)_{(Y)}\, \Pi&=
\prod_{k=1}^{n+1} \frac{1-qv_k/tu_i}{1-qv_k/u_i}
\, H(u_1,\ldots,u_i/q,\ldots,u_{n+1};v),
\end{split}
\end{equation*}
and (ii) follows from (3).
\end{proof}

For any partition $\la=(\la_1,\ldots,\la_{n})$ and 
$1 \le k \le n$, we write 
$\la_{(k)}=(\la_1,\ldots,\la_k-1,\ldots,\la_{n})$.
\begin{Lemma}
Let $X=(x_1,\ldots,x_{N})$ and $z$ be $N+1$ indeterminates and 
$\la=(\la_1,\ldots,\la_{n})$ an arbitrary partition with length
$n \le N$. The coefficient of  $z$ in
$Q_{\la}(X,z;q,t)$ is given by
\[\frac{1-t}{1-q}\, \sum_{k=1}^n c_k(\la) \,
Q_{\la_{(k)}}(X;q,t),\]
with $c_k(\la)=0$ if $\la_{(k)}$ is not a partition, and otherwise
\[c_k(\la)=\prod_{i=1}^{k-1}\frac{1-q^{\la_i-\la_k+1}t^{k-i-1}}
{1-q^{\la_i-\la_k+1}t^{k-i}} \,
\frac{1-q^{\la_i-\la_k}t^{k-i+1}}
{1-q^{\la_i-\la_k}t^{k-i}}.\]
\end{Lemma}

\begin{proof}
We have
\[Q_{\la}(X,z)=\sum_{\mu}Q_{\la/\mu}(X)\,Q_{\mu}(z),\]
where $Q_{\la/\mu}$ denotes the skew Macdonald 
function \cite[p.~345]{Ma}.
But $Q_{\mu}(z)$ is non zero only if $\mu$ is a 
row partition $(k)$, in which case $Q_{(k)}(z)=g_k(z)=a_kz^k$ for 
some constant $a_k$. In particular $a_1=(1-t)/(1-q)$.
Thus the coefficient of $z$ in $Q_{\la}(X,z)$ is
$(1-t)/(1-q) \,Q_{\la/(1)}(X)$.
Now using Macdonald's notations, the skew function 
$Q_{\la/(1)}$ is defined by
\[Q_{\la/(1)}=\sum_{\nu}\psi^{'}_{\la/\nu}\,Q_{\nu},\]
with $\psi^{'}_{\la/\nu}$ given by
\[e_1\,P_{\nu}=\sum_{\la}\psi^{'}_{\la/\nu}\,P_{\la}.\]
It is the simplest case of 
the Pieri formula for Macdonald polynomials. The quantities 
$\psi^{'}_{\la/\nu}$ have been computed in \cite[p.~336]{Ma}. 
The partition $\nu$ writes necessarily as
$\nu=\la_{(k)}=(\la_1,\ldots,\la_k-1,\ldots,\la_{n})$ for 
$1\le k \le n$, and we have 
$\psi^{'}_{\la/\la_{(k)}}=c_k(\la)$.
\end{proof}
 
\begin{proof}[Proof of Theorem 1]
We start from the Pieri formula proved in \cite[p.~340, 
(6.24), (ii)]{Ma}. Being given a partition 
$\mu=(\mu_1,\ldots,\mu_n)$ with length 
$n$, this Pieri formula reads
\[ Q_{(r)} \: Q_{\mu}=\sum_{\kappa \supset\mu} 
\psi_{\kappa/\mu}\,Q_{\kappa},\]
where the skew diagram $\kappa-\mu$ is a horizontal 
$r$ - strip, i.e. has at most one square in each column.
Since $\kappa-\mu$ is a horizontal strip, the length of 
$\kappa$ is at most equal to $n+1$, and we can write
\begin{equation*}
Q_{(r)} \:
Q_{(\mu_1,\ldots,\mu_n)}= \sum_{\ta\in\mathsf{N}^n}
\psi_{\ta_1,\ldots,\ta_n} (\mu)\:
Q_{(\mu_1+\ta_1,\ldots,\mu_n+\ta_n,r-|\ta|)}.
\end{equation*}

Now let  $\la=(\la_1,...,\la_{n+1})$ be a partition 
with length $n+1$. Expanding $Q_{(\la_{n+1})}\,
Q_{(\la_1,\ldots,\la_n)}$ and using induction 
on $\la_{n+1}$, we see that the Pieri 
formula can be inverted as
\begin{equation*}
Q_{(\la_1,\ldots,\la_{n+1})}= \sum_{\ta\in\mathsf{N}^n}
F_{\ta_1,\ldots,\ta_n} (\la)\:
Q_{(\la_{n+1}-|\ta|)} \:
Q_{(\la_1+\ta_1,\ldots,\la_n+\ta_n)},
\end{equation*}
with $F_{\ta} (\la)$ to be determined.
Since $\psi_{0,\ldots,0} (\mu)=1$, one has $F_{0,\ldots,0} (\la)=1$.

We can write this expansion at $(X,z)$, apply Lemma 2, and 
identify the coefficient of $z$ on each side, up to 
$(1-t)/(1-q)$. The 
coefficient on the left-hand side is 
\begin{multline*}
\sum_{k=1}^{n+1} c_k(\la) \,
Q_{\la_{(k)}}= \sum_{\ta\in\mathsf{N}^n} \Big(
\sum_{k=1}^{n} c_k(\la) \,
F_{\ta_1,\ldots,\ta_n} (\la_{(k)}) \,
Q_{(\la_{n+1}-|\ta|)} \:
Q_{(\la_1+\ta_1,\ldots,\la_k+\ta_k-1,\ldots,\la_n+\ta_n)}\\ +
c_{n+1}(\la) \, 
F_{\ta_1,\ldots,\ta_n} (\la_{(n+1)}) \,
Q_{(\la_{n+1}-|\ta|-1)} \:
Q_{(\la_1+\ta_1,\ldots,\la_n+\ta_n)}\Big).
\end{multline*}
The coefficient on the right-hand side is 
\begin{multline*}
\sum_{\ta\in\mathsf{N}^n}
F_{\ta_1,\ldots,\ta_n} (\la)
\Big(Q_{(\la_{n+1}-|\ta|-1)} \:
Q_{(\la_1+\ta_1,\ldots,\la_n+\ta_n)}\\ +
\sum_{k=1}^{n} c_k((\la_1+\ta_1,\ldots,\la_n+\ta_n)) \,
Q_{(\la_{n+1}-|\ta|)} \:
Q_{(\la_1+\ta_1,\ldots,\la_k+\ta_k-1,\ldots,\la_n+\ta_n)}\Big).
\end{multline*}
Now the products $Q_{(r)}Q_{\mu}$ 
form a basis of $\mathsf{Sym}$.
Identifying coefficients of 
$Q_{(\la_{n+1}-|\ta|-1)}$
$Q_{(\la_1+\ta_1,\ldots,\la_n+\ta_n)}$ on both sides, 
we obtain
\begin{multline}
\sum_{i=1}^n c_i(\la)
\, F_{\ta_1,\ldots,\ta_i+1,\ldots,\ta_n} (\la_{(i)})
+c_{n+1}(\la) F_{\ta_1,\ldots,\ta_n} 
(\la_{(n+1)})=\\
F_{\ta_1,\ldots,\ta_n}(\la) + \sum_{i=1}^n
c_i((\la_1+\ta_1,\ldots,\la_i+\ta_i+1,\ldots,\la_n+\ta_n)) \,
F_{\ta_1,\ldots,\ta_i+1,\ldots,\ta_n} (\la).
\end{multline}
Remark that in order to cover all the cases, this equation 
must also be considered when one $\ta_k$ is equal to $-1$.
Then the corresponding $F_{\ta} (\la)$ are of course set to zero. 

An essential fact is that equation (4) can be used to define 
$F_{\ta_1,\ldots,\ta_n} (\la)$ by induction, starting from 
$F_{0,\ldots,0} (\la)=1$. Indeed it can be written as
\begin{multline*}
c_n((\la_1+\ta_1,\ldots,\la_n+\ta_n)) \,
F_{\ta_1,\ldots,\ta_n} (\la)=-F_{\ta_1,\ldots,\ta_n-1}(\la)\\
-\sum_{i=1}^{n-1}
c_i((\la_1+\ta_1,\ldots,\la_i+\ta_i+1,\ldots,\la_n+\ta_n-1)) \,
F_{\ta_1,\ldots,\ta_i+1,\ldots,\ta_n-1} (\la)\\
+ \sum_{i=1}^n c_i(\la)
\, F_{\ta_1,\ldots,\ta_i+1,\ldots,\ta_n-1} (\la_{(i)})
+c_{n+1}(\la) F_{\ta_1,\ldots,\ta_n-1} (\la_{(n+1)}).
\end{multline*}

This defines $F_{\ta} (\la)$ through an induction on $|\la|$ 
and $\ta_n$. Indeed each $F_{\kappa} (\mu)$ which contributes to 
the right-hand side has either $|\mu|=|\la|-1$, or 
$\kappa_n=\ta_n-1$. And for $\ta_n=0$ the previous relation 
writes
\[c_n((\la_1+\ta_1,\ldots,\la_{n-1}+\ta_{n-1},\la_n)) \,
F_{\ta_1,\ldots,\ta_{n-1},0} (\la)=
c_n(\la)\, F_{\ta_1,\ldots,\ta_{n-1},0} (\la_{(n)}),\]
which still allows induction on $|\la|$.

Moreover since $c_k(\la)$ depends only on the quantities $q^{\la_i-\la_j}$, 
and not on $q^{\la_{n+1}}$, the inductive definition implies that 
$F_{\ta_1,\ldots,\ta_n} (\la)$ depends only on 
$q^{\la_{i}-\la_{n+1}}, 1\le i\le n$. So we may write
\[F_{\ta_1,\ldots,\ta_n}(\la_1,\ldots,\la_{n+1})=
F_{\ta_1,\ldots,\ta_n}(u_1,\ldots,u_n),\]
with $u_i=q^{\la_i-\la_{n+1}}t^{n-i}$, $1\le i \le n$. 
Setting $v_i=q^{\ta_i}u_i$, 
relation (4) may be written
\begin{multline}
\sum_{i=1}^n \prod_{k=1}^{i-1}
\frac{1-qu_k/tu_i}{1-qu_k/u_i} \,
\frac{1-tu_k/u_i}{1-u_k/u_i}
\, F_{\ta_1,\ldots,\ta_i+1,\ldots,\ta_n}
(u_1,\ldots,u_i/q,\ldots,u_n)\\
+\prod_{k=1}^{n} \frac{1-qu_k}{1-qtu_k} \,
\frac{1-t^2u_k}{1-tu_k} \,
F_{\ta_1,\ldots,\ta_n} (qu)=\\ 
F_{\ta_1,\ldots,\ta_n} (u) + \sum_{i=1}^n \prod_{k=1}^{i-1}
\frac{1-v_k/tv_i}{1-v_k/v_i} \,
\frac{1-tv_k/qv_i}{1-v_k/qv_i} \,
F_{\ta_1,\ldots,\ta_i+1,\ldots,\ta_n}(u).
\end{multline}

By the inductive hypothesis, the property being true for $\ta=0$, 
the proof will be finished if we show that 
$C^{(q,t)}_{\ta_1,\ldots,\ta_n} (u_1,\ldots,u_n)$ satisfies 
the previous equation. Using
\begin{multline*}
C^{(q,t)}_{\ta_1,\ldots,\ta_n} (u_1,\ldots,u_n)=\\
\prod_{1\le i< j\le n+1}
\frac{(qu_i/tu_j;q)_{\ta_i}}{(qu_i/u_j;q)_{\ta_i}}
\prod_{1\le i\le j\le n}
\frac{(tu_i/v_j;q)_{\ta_i}}{(u_i/v_j;q)_{\ta_i}}\:
\lim_{a\rightarrow-1}\frac{1}{1+a} H(u,1/t;v,0),
\end{multline*}
after an easy computation, it is sufficient to prove
\begin{multline*}
\sum_{i=1}^n \frac{1-u_i}{1-tu_i} \,
\prod_{\begin{subarray}{c}k=1\\ k\neq i\end{subarray}}^{n}
\frac{u_i/t-u_k}{u_i-u_k} \,
\prod_{k=1}^n
\frac{tu_i-qv_k}{u_i-qv_k}
\, H(u_1,\ldots,u_i/q,\ldots,u_n,1/t;v,0)\\
+\prod_{k=1}^{n} \frac{1-t^2u_k}{1-tu_k} \,
\frac{1-qv_k}{1-qtv_k}\, 
H(qu,1/t;qv,0)= H(u,1/t;v,0) +\\
\sum_{i=1}^n \frac{1-qv_i}{1-qtv_i} \,
\prod_{\begin{subarray}{c}k=1\\ k\neq i\end{subarray}}^{n}
\frac{v_i-v_k/t}{v_i-v_k} \,
\prod_{k=1}^n
\frac{tu_k-qv_i}{u_k-qv_i} \,
H(u,1/t;v_1,\ldots,qv_i,\ldots,v_n,0).
\end{multline*}
Since $H(u,u_{n+1}/q;v,0)=H(qu,u_{n+1};qv,0)$ 
we conclude by applying Lemma 1 (ii), with $u_{n+1}=1/t$ 
and $v_{n+1}=0$. 
\end{proof}

\noindent \textit{Remark :} Starting from Lemma 1 (i), the 
same argument applied in reverse order shows 
that $F_{\ta_1,\ldots,\ta_n}=C^{(q,t)}_{\ta_1,\ldots,\ta_n}$ 
also satisfies the recurrence relation
\begin{multline*}
\sum_{k=1}^n \prod_{i=k+1}^{n}
\frac{1-v_k/tv_i}{1-v_k/v_i} \,
\frac{1-tv_k/qv_i}{1-v_k/qv_i}
\, F_{\ta_1,\ldots,\ta_k-1,\ldots,\ta_n}(u)
+ F_{\ta_1,\ldots,\ta_n} (u)= 
F_{\ta_1,\ldots,\ta_n} (u/q) +\\
\sum_{k=1}^n 
\frac{1-qu_k}{1-qtu_k} \,
\frac{1-t^2u_k}{1-tu_k} \,
\prod_{i=k+1}^{n}
\frac{1-qu_k/tu_i}{1-qu_k/u_i} \,
\frac{1-tu_k/u_i}{1-u_k/u_i} \,
F_{\ta_1,\ldots,\ta_k-1,\ldots,\ta_n}(u_1,\ldots,qu_k,\ldots,u_n).
\end{multline*}
However a direct proof of this second relation seems much more difficult
than for (5).

\end{document}